\documentclass[12pt,twoside]{article}
\usepackage{amssymb,amsmath,amsthm,latexsym,graphics,graphicx,enumerate,cite}
\usepackage{amscd}
\usepackage[cp1251]{inputenc} \usepackage[english]{babel}

\newcommand{\wt}[1]{\widetilde{#1}}
\newcommand{\mc}[1]{\mathcal{#1}}
\newcommand{\mb}[1]{\mathbb{#1}}
\newcommand{\ov}[1]{\overline{#1}}
\newcommand{\vp}{\varphi}

\newcommand{\RDf}{\mathrm{Rep}\,(\widetilde{D}_4,f)}
\newcommand{\RPf}{\mathrm{Rep}\,{\mathcal{P}_{4,f}}}

\DeclareMathOperator{\Imm}{Im}

\begin{document}

\makeatletter
\renewcommand{\@biblabel}[1]{#1.}
\makeatother

\begin{center}
\Large \bf
On Regular Locally Scalar Representations\\
of Graph $\widetilde{D}_4$ in Hilbert Spaces
\end{center}

\begin{center}
\large\bfseries S.~A.~Kruglyak~$^\dag$, L.~A.~Nazarova, \fbox{A.~V.~Roiter}.
\end{center}

\noindent
$^\dag$~Institute of Mathematics of National Academy of Sciences of Ukraine,\\
Tereshchenkovska str., 3, Kiev, Ukraine, ind. 01601\\
E-mail: red@imath.kiev.ua

\bigskip

Representations of quivers corresponding to extended Dynkin graphs are
described up to equivalence in \cite{Naz1}. Locally scalar
representations of graphs in the category of Hilbert spaces were introduced
in \cite{KruRoi1}, and such representations are naturally classified
up to unitary equivalense.

Representations of $*$-algebras generated by linearly related orthogonal
projections are studied in
\cite{Bes1,OstSam1,GalKru1,Gal1,GalMur1,KruRabSam1,KruKir1,Str1} and
others. The connection between locally scalar representations of several
graphs (trees, which include also Dynkin graphs) and representations of
such $*$-algebras is stated in \cite{KruPopSam1}, and we further use this
connection.

The present paper is dedicated to the classification of indecomposable
regular (see \cite{RedRoi1}) locally scalar representations of the graph
$\wt{D}_4$ (for $\wt{D}_4$ those are indecomposable locally scalar
representations in the dimension $(2;1,1,1,1)$). The answer obtained for
the corresponding $*$-algebra in \cite{GalKru1,Gal1}, in our opinion,
cannot be satisfactory and definitive. We will obtain explicit formulas,
expressing matrix elements of a representation by a character (see
\cite{KruRoi1}) of a locally scalar representation and two ``free'' real
parameters.

\bigskip

{\bf 1.}
Let $\mc{H}$ be the category of finite-dimensional Hilbert spaces, and
$\wt{D}_4$ be an extended Dynkin graph

\centerline{\includegraphics[scale=0.85]{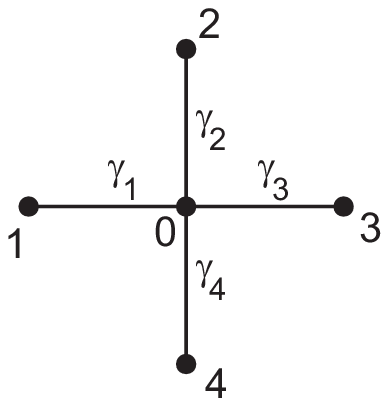}}

Remind (see \cite{KruRoi1}) that a representation $\Pi$ of the graph
$\wt{D}_4$ associates a space $H_i\in\mc{H}$ to each vertex $i$ ($i=\ov{0,4}$),
and a pair of interadjoint linear operators
$\Pi(\gamma_i)=\{\Gamma_{i0};\Gamma_{0i}\}$ to each edge $\gamma_i$, where
$\Gamma_{0i}: H_i\to H_0$, $\Gamma_{i0}=\Gamma_{0i}^*$.

A morphism $C: \Pi\to\wt{\Pi}$ is a family $C=\{C_i\}_{i=\ov{0,4}}$ of
operators $C_i: \Pi(i)\to\wt{\Pi}(i)$ such that the diagrams
$$
\begin{CD}
H_i @>\Gamma_{ji}>> H_j\\
@V C_i VV @V C_j VV\\
\wt{H}_i @>\wt{\Gamma}_{ji}>> \wt{H}_j
\end{CD}
$$
\noindent are commutative, i.~e. $C_j\Gamma_{ji}=\wt{\Gamma}_{ji}C_i$.

Let $M_i$ be the set of vertices connected with a vertex $i$ by an edge,
$A_i=\sum_{j\in M_i} \Gamma_{ij}\Gamma_{ji}$. Representation $\Pi$ is
called locally scalar \cite{KruRoi1} if all operators $A_i$ are scalar;
$A_i=\alpha_iI_{H_i}$, where $I_{H_i}$ is identity operator in a space
$H_i$. Since $A_i$ is a positive operator, $\alpha_i\geq 0$. A vector
$\{\dim \Pi(i)\}$ is a dimension of a finite-dimensional representation
$\Pi_i$; if $A_i=f(i)I_{H_i}$ then $\{f(i)\}$ is called a character of
locally scalar representation $\Pi$.

Further we will denote as $\RDf$ the category of finite-dimensional locally
scalar representations of the graph $\wt{D}_4$ in $\mc{H}$ with given
character $f$.

We will assume that $\alpha_i=f(i)>0,\; i=\ov{0,4}$ and the character is
normalized: $f(0)=\alpha_0=1$.

On the other hand, consider the following $*$-algebra over the field
$\mb{C}$:
$$
{\mc{P}_{4,f}} = \mb{C}\langle
p_1,p_2,p_3,p_4 \,|\, p_i=p_i^*=p_i^2,\;
\sum_{i=1}^4 \alpha_ip_i=e
\rangle,
$$
\noindent where $\alpha_i=f(i)$, $e$ is the identity of the algebra, and
the category $\RPf$ of finite-dimensional $*$-representation of the algebra
$\mc{P}_{4,f}$.

Let $\Pi\in\RDf,\; \Pi(i)=H_i,\;
\pi(\gamma_i)=\{\Gamma_{i,0};\Gamma_{0,i}\}$. Let us construct a
representation $\pi$ of the algebra $\mc{P}_{4,f}$ by the following way:
$\pi(p_i)=\frac{1}{\alpha_i}\Gamma_{0,i}\cdot\Gamma_{i,0}=P_i$. If $C:
\Pi\to\wt{\Pi}$ is a morphism in $\RDf$ then $C_0: \pi\to\wt{\pi}$ is a
morphism in the category $\RPf$ ($C_0$ is the operator interlacing
representations $\pi$ and $\wt{\pi}$). Define a functor $$ \Phi:
\RDf\to\RPf $$ \noindent putting $\Phi(\Pi)=\pi,\; \Phi(C)=C_0$. Clearly,
the functor $\Phi$ is the equivalence of categories.

Let $\pi\in\RPf$ be a representation in the space $H_0$. Set $H_i = \Imm
P_i,\; i=\ov{1,4};\; \Gamma_{0,i}: H_i\to H_0$ is the natural injection of
the space $H_i$ into $H_0$, then, putting
$\Pi(\gamma_i)=\{\Gamma_{0,i};\Gamma_{0,i}^*\}$, we obtain a representation
from $\RDf$. If $C_0:\pi\to\wt{\pi}$ set $C_i=C_0\big\vert_{H_i}$. If
$\Phi^{(-1)}(\pi)=\Pi,\; \Phi^{(-1)}(C_0)=\{C_i\}_{i=\ov{0,4}}$ then
$\Phi\Phi^{(-1)} \sim I_{\RPf},\; \Phi^{(-1)}\Phi \sim I_{\RDf}$.

\bigskip

{\bf 2.} Consider representations of the $*$-algebra $\mc{P}_{4,f}$ for
$$
0<\alpha_1\leq\alpha_2\leq\alpha_3\leq\alpha_4< 1,\;
\sum_{i=1}^4\alpha_i=2.
$$
\noindent (in the other cases representations of the $*$-algebras
$\mc{P}_{4,f}$ are reduced to the simplest representations by the Coxeter
functors \cite{Kru1}).

Let us make the substitution of generators in the algebra $\mc{P}_{4,f}$
\cite{GalKru1} :
\begin{equation}\label{eq01}
\begin{array}{l}
x=\alpha_2p_2+\alpha_3p_3-\frac{1}{2}\beta_1e,\\
y=\alpha_1p_1+\alpha_3p_3-\frac{1}{2}\beta_2e,\\
z=\alpha_1p_1+\alpha_2p_2-\frac{1}{2}\beta_3e,
\end{array}
\text{ where }
\begin{array}{l}
\beta_1=(2-\alpha_1+\alpha_2+\alpha_3-\alpha_4)/2,\\
\beta_2=(2+\alpha_1-\alpha_2+\alpha_3-\alpha_4)/2,\\
\beta_3=(2+\alpha_1+\alpha_2-\alpha_3-\alpha_4)/2.
\end{array}
\end{equation}
Denote also
$$
\begin{array}{l}
\gamma_1=(\alpha_1^2-\alpha_2^2-\alpha_3^2+\alpha_4^2)/4,\\
\gamma_2=(-\alpha_1^2+\alpha_2^2-\alpha_3^2+\alpha_4^2)/4,\\
\gamma_3=(-\alpha_1^2-\alpha_2^2+\alpha_3^2+\alpha_4^2)/4.
\end{array}
$$
It is easy to check that $\gamma_1\leq\gamma_2\leq\gamma_3$ and
$0\leq\gamma_2$.

The new generators $x,\,y,\,z$ satisfy the system of relations
\begin{equation}\label{eq02}
\begin{array}{l}
\{y,z\}=\gamma_1e,\\
\{z,x\}=\gamma_2e,\\
\{x,y\}=\gamma_3e,
\end{array}\quad
(x+y+z)^2=\alpha_4^2e.
\end{equation}
The equalities \eqref{eq01} imply
\begin{equation}\label{eq03}
\begin{array}{l}
p_1=\cfrac{-x+y+z}{2\alpha_1}+\frac{1}{2}e,\\
p_2=\cfrac{x-y+z}{2\alpha_2}+\frac{1}{2}e,\\
p_3=\cfrac{x+y-z}{2\alpha_3}+\frac{1}{2}e,\\
p_4=\cfrac{-x-y-z}{2\alpha_4}+\frac{1}{2}e.
\end{array}
\end{equation}

\bigskip

{\bf 3.} Let $\gamma_3=0$, then $0\leq\gamma_2\leq\gamma_3$ implies
$\gamma_2=0$, hence $\alpha_1=\alpha_2=\alpha_3=\alpha_4=\frac{1}{2}$ (this
case is considered in \cite{OstSam1}), and so $\gamma_1=0$.

In this case the system of relation \eqref{eq02} has a form
\begin{equation}\label{eq04}
\begin{array}{l}
\{y,z\}=0,\\
\{z,x\}=0,\\
\{x,y\}=0,
\end{array}\quad
x^2+y^2+z^2=\frac{1}{4}e.
\end{equation}
Let $\pi$ be an indecomposable two-dimensional representation of the algebra
$\mathcal{P}_{4,f}$ and
$\pi(x)=X,\; \pi(y)=Y,\; \pi(z)=Z$.
\begin{itemize}
\item[a)] Let $Z=0$. The matrix $X$ can be diagonalized as a matrix of
self-adjoint operator. The relations of anticommutation imply that the
triple $X,\,Y,\,Z=0$ is indecomposable only in the case when the
diagonalized matrix $X$ equals
$\begin{bmatrix} -\lambda& 0\\ 0& \lambda \end{bmatrix}$.
Then $Y=\begin{bmatrix} 0& y_{12}\\ \ov{y}_{12}& 0 \end{bmatrix}$,
$y_{12}\neq 0$, and the element $y_{12}$ can be made positive by the
admissible tranformations. Therefore, we obtain the case
$$
X=\lambda\begin{bmatrix}-1&0\\0&1\end{bmatrix},\quad
Y=\mu\begin{bmatrix}0&1\\1&0\end{bmatrix},\quad
Z=\begin{bmatrix}0&0\\0&0\end{bmatrix};\quad
\lambda >0,\, \mu >0.
$$
\item[b)] Let $Z\neq 0$ and $X=0$. Then the matrix $Y$ can be diagonalized:
$Y=\begin{bmatrix}-\mu&0\\0&\mu\end{bmatrix},\; \mu >0$ (in the other cases
the triple of matrices turn out to be decomposable). Then from $\{y,z\}=0$
we obtain
$Z=\begin{bmatrix}0&z_{12}\\\ov{z}_{12}&0\end{bmatrix}$ and one can reduce
the matrix $Z$ (does not changing the $Y$) by the admissible transformations
to the form
$Z=\begin{bmatrix}0&i\nu\\\-i\nu&0\end{bmatrix},\; \nu > 0$.
The triple of matrices $X,\,Y,\,Z$ is reduced by means
of a unitary matrix $U=\begin{bmatrix}\sqrt{2}/2&\sqrt{2}/2\\
-\sqrt{2}/2&\sqrt{2}/2\end{bmatrix}$ by the transormation
$UXU^*,\,UYU^*,\,UZU^*$ to the form
$$
X=\begin{bmatrix}0&0\\0&0\end{bmatrix},\quad
Y=\mu\begin{bmatrix}0&1\\1&0\end{bmatrix},\quad
Z=\nu\begin{bmatrix}0&i\\-i&0\end{bmatrix};\quad
\mu >0,\, \nu >0.
$$
\item[c)] Let $X\neq 0,\; Z\neq 0,\; Y=0$. In this case the matrices
$X,\,Y,\,Z$ can be reduced to the form
$$
X=\lambda\begin{bmatrix}-1&0\\0&1\end{bmatrix},\quad
Y=\begin{bmatrix}0&0\\0&0\end{bmatrix},\quad
Z=\nu\begin{bmatrix}0&i\\-i&0\end{bmatrix};\quad
\lambda >0,\, \nu >0.
$$
\item[d)] $X\neq 0,\; Y\neq 0,\; Z\neq 0$. In this case the matrices
$X,\,Y,\,Z$ can be reduced to the form
$$
X=\lambda\begin{bmatrix}-1&0\\0&1\end{bmatrix},\quad
Y=\mu\begin{bmatrix}0&1\\1&0\end{bmatrix},\quad
Z=\nu\begin{bmatrix}0&i\\-i&0\end{bmatrix};\quad
\lambda >0,\, \mu >0,\, \nu\in\mb{R},\, \nu\neq 0.
$$
\end{itemize}
Thus, unifying these cases, we may consider that
$$
\begin{array}{l}
X=\lambda\begin{bmatrix}-1&0\\0&1\end{bmatrix},\quad
Y=\mu\begin{bmatrix}0&1\\1&0\end{bmatrix},\quad
Z=\nu\begin{bmatrix}0&i\\-i&0\end{bmatrix},\\
\text{where}\; \lambda^2+\mu^2+\nu^2=\frac{1}{4}\;
\text{(follows from \eqref{eq04})},\;\text{and}\\
	\begin{array}{l}
	\qquad\text{either}\; \lambda >0,\,\mu >0,\,\nu\in\mb{R};\\
	\qquad\text{either}\; \lambda =0,\,\mu >0,\,\nu >0;\\
	\qquad\text{or}\; \lambda >0,\,\mu =0,\,\nu >0.
	\end{array}
\end{array}
$$
Formulas \eqref{eq03} imply
\begin{eqnarray*}
P_1=\begin{bmatrix}
\frac{1}{2}-\lambda & \mu+\nu i\\
\mu-\nu i & \frac{1}{2}+\lambda
\end{bmatrix},&
P_2=\begin{bmatrix}
\frac{1}{2}+\lambda & -\mu+\nu i\\
-\mu-\nu i & \frac{1}{2}-\lambda
\end{bmatrix},\\
P_3=\begin{bmatrix}
\frac{1}{2}+\lambda & \mu-\nu i\\
\mu+\nu i & \frac{1}{2}-\lambda
\end{bmatrix},&
P_4=\begin{bmatrix}
\frac{1}{2}-\lambda & -\mu-\nu i\\
-\mu+\nu i & \frac{1}{2}+\lambda
\end{bmatrix}.
\end{eqnarray*}
If $\mu+\nu i=\sqrt{\mu^2+\nu^2}\,e^{i\vp}=\sqrt{\frac{1}{4}-\lambda^2}\,e^{i\vp}$
then, passing on to the unitary equivalent representation  by means of the
matrix $U=\begin{bmatrix}e^{-ivp}&0\\0&1\end{bmatrix}$, we may consider that
\begin{eqnarray}\label{eq05}
P_1=\begin{bmatrix}
\frac{1}{2}-\lambda & \sqrt{\frac{1}{4}-\lambda^2}\\
\sqrt{\frac{1}{4}-\lambda^2} & \frac{1}{2}+\lambda
\end{bmatrix},&
P_2=\begin{bmatrix}
\frac{1}{2}+\lambda & e^{i\chi}\sqrt{\frac{1}{4}-\lambda^2}\\
e^{-i\chi}\sqrt{\frac{1}{4}-\lambda^2}& \frac{1}{2}-\lambda
\end{bmatrix},\nonumber\\
P_3=\begin{bmatrix}
\frac{1}{2}+\lambda & -e^{-i\chi}\sqrt{\frac{1}{4}-\lambda^2}\\
-e^{-i\chi}\sqrt{\frac{1}{4}-\lambda^2} & \frac{1}{2}-\lambda
\end{bmatrix},&
P_4=\begin{bmatrix}
\frac{1}{2}-\lambda & -\sqrt{\frac{1}{4}-\lambda^2}\\
-\sqrt{\frac{1}{4}-\lambda^2} & \frac{1}{2}+\lambda
\end{bmatrix},
\end{eqnarray}
$$
\begin{array}{c}
0\leq\lambda < 1/2,\; \text{and}\\
\qquad 0<\chi <\pi/2\; \text{when}\; \lambda = 0,\\
\qquad -\pi/2<\chi \leq\pi/2\; \text{when}\; 0<\lambda <1/2.
\end{array}
$$

\bigskip

{\bf 4.} Let $\gamma_3\neq 0$, then $\{x,y\}=\gamma_3e$ implies
$X\neq 0$ and $Y\neq 0$. Moreover, $X$ and $Y$ has no zero eigenvalues.
Indeed, let martix $X$ has the form after digonalization:
$X=\begin{bmatrix}\lambda_1&0\\0&\lambda_2\end{bmatrix}$. In the matrices
$Y$ and $Z$ either $y_{12}\neq 0$ or $z_{12}\neq 0$ (or else the triple of
matrices is decomposable). $\{x,y\}=\gamma_3e,\; \{x,z\}=\gamma_2e$ imply
$(\lambda_1+\lambda_2)z_{12}=0$ and $(\lambda_1+\lambda_2)y_{12}=0$.
Therefore we can conclude that $-\lambda_1=\lambda_2=\lambda\neq 0,\;
\lambda >0$. The same reasoning is useful also for $Y$.

Let $X=\begin{bmatrix}-\lambda&0\\0&\lambda\end{bmatrix}$.
$\{x,y\}=\gamma_3e$ implies $y_{11}=-\frac{\gamma_3}{2\lambda},\;
y_{22}=\frac{\gamma_3}{2\lambda}$;
$\{x,z\}=\gamma_2e$ implies $z_{11}=-\frac{\gamma_2}{2\lambda},\;
z_{22}=\frac{\gamma_2}{2\lambda}$:
$$
X=\lambda\begin{bmatrix}-1&0\\0&1\end{bmatrix},\quad
Y=\frac{1}{2\lambda}\begin{bmatrix}
-\gamma_3&y_{12}\\\ov{y}_{12}&\gamma_3\end{bmatrix},\quad
Z=\frac{1}{2\lambda}\begin{bmatrix}
-\gamma_2&z_{12}\\\ov{z}_{12}&\gamma_2\end{bmatrix}.
$$
$\{y,z\}=\gamma_1e$ implies
$$
\frac{1}{4\lambda^2}\begin{bmatrix}
2\gamma_2\gamma_3+y_{12}\ov{z}_{12}+\ov{y}_{12}z_{12}&0\\
0&2\gamma_2\gamma_3+y_{12}\ov{z}_{12}+\ov{y}_{12}z_{12}
\end{bmatrix}=\gamma_1I,
$$
\noindent hence
\begin{equation}\label{eq06}
y_{12}\ov{z}_{12}+\ov{y}_{12}z_{12}=4\gamma_1\lambda^2-2\gamma_2\gamma_3.
\end{equation}
Let us turn to the unitary equivalent representation by means of the
unitary matrix $U=\begin{bmatrix}e^{i\vp}&0\\0&0\end{bmatrix}$ so that
$y_{12}+z_{12}=r_1$ would be real positive; at that $-y_{12}+z_{12}$
remains to be complex in general, $-y_{12}+z_{12}=r_2e^{i\chi}$.

Then \eqref{eq02} implies
$$
(X+Y+Z)^2=\frac{1}{4\lambda^2}
\begin{bmatrix}
(2\lambda^2+\gamma_2+\gamma_3)^2+(y_{12}+z_{12})^2&0\\
0&(2\lambda^2+\gamma_2+\gamma_3)^2+(y_{12}+z_{12})^2.
\end{bmatrix}=\alpha_4^2I
$$
\noindent and
$$
r_1^2=(y_{12}+z_{12})^2=
4\alpha_4^2\lambda^2-(2\lambda^2+\gamma_2+\gamma_3)^2,
$$
\noindent from which it easy to obtain
\begin{equation}\label{eq07}
\begin{array}{ll}
r_1&=
\sqrt{-4\lambda^4+2(\alpha_1^2+\alpha_4^2)\lambda^2-
\frac{1}{4}(\alpha_4^2-\alpha_1^2)^2}=\\
&=\sqrt{-4\left(\lambda^2-\frac{(\alpha_4^2-\alpha_1^2)^2}{4}\right)
\left(\lambda^2-\frac{(\alpha_4^2+\alpha_1^2)^2}{4}\right)},
\end{array}
\end{equation}
$$
\frac{\alpha_4-\alpha_1}{2}\leq\lambda\leq\frac{\alpha_4+\alpha_1}{2}.
$$
$$
\begin{cases}
y_{12}+z_{12}=r_1,\\
-y_{12}+z_{12}=r_2e^{i\chi},
\end{cases}
\;\text{therefore}\;
\begin{cases}
y_{12}=(r_1-r_2e^{i\chi})/2,\\
z_{12}=(r_1+r_2e^{i\chi})/2,
\end{cases}
$$
\noindent $y_{12}\ov{z}_{12}+\ov{y}_{12}z_{12}=(r_1^2-r_2^2)/2=
4\gamma_1\lambda^2-2\gamma_2\gamma_3$ (the last equality follows from
\eqref{eq06}).

Thus,
$$
r_2^2=r_1^2+4\gamma_2\gamma_3-8\gamma_1\lambda^2=
-4\lambda^4+2(\alpha_2^2+\alpha_3^2)\lambda^2-
\frac{1}{4}(\alpha_3^2-\alpha_2^2)^2\geq 0
$$
\noindent and
$$
\begin{array}{ll}
-y_{12}+z_{12}=r_2e^{i\chi}&=
e^{i\chi}\sqrt{-4\lambda^4+2(\alpha_2^2+\alpha_3^2)\lambda^2-
\frac{1}{4}(\alpha_3^2-\alpha_2^2)^2}=\\
&=e^{i\chi}\sqrt{-4\left(\lambda^2-\frac{(\alpha_3^2-\alpha_2^2)^2}{4}\right)
\left(\lambda^2-\frac{(\alpha_3^2+\alpha_2^2)^2}{4}\right)},
\end{array}
$$
$$
\frac{\alpha_3-\alpha_2}{2}\leq\lambda\leq\frac{\alpha_3+\alpha_2}{2}.
$$
Now we can directly pass on to the determining of the operators of
the representation $P_1,\,P_2,\,P_3,\,P_4$:
\begin{multline}
P_1=\cfrac{\!-\!X\!+\!Y\!+\!Z}{2\alpha_1}+\frac{1}{2}I = \\
=\cfrac{1}{4\alpha_1\lambda}
\begin{bmatrix}
2\lambda^2\!+\!2\alpha_1\lambda\!-\!\frac{1}{2}(\alpha_4^2\!-\!\alpha_1^2)&
\sqrt{\!-\!4\lambda^4\!+\!2(\alpha_1^2\!+\!\alpha_4^2)\lambda^2\!-\!
\frac{1}{4}(\alpha_4^2\!-\!\alpha_1^2)^2}\\
\sqrt{\!-\!4\lambda^4\!+\!2(\alpha_1^2\!+\!\alpha_4^2)\lambda^2\!-\!
\frac{1}{4}(\alpha_4^2\!-\!\alpha_1^2)^2}&
\!-\!2\lambda^2\!+\!2\alpha_1\lambda\!+\!\frac{1}{2}(\alpha_4^2\!-\!\alpha_1^2)
\end{bmatrix}=\\
=\cfrac{1}{4\alpha_1\lambda}
\begin{bmatrix}
2\left(\lambda\!-\!\frac{1}{2}(\alpha_4\!-\!\alpha_1)\right)
\left(\lambda\!+\!\frac{1}{2}(\alpha_4\!+\!\alpha_1)\right)&
\sqrt{\!-\!4\left(\lambda^2\!-\!\frac{1}{4}(\alpha_4\!-\!\alpha_1)^2\right)
\left(\lambda^2\!-\!\frac{1}{4}(\alpha_4\!+\!\alpha_1)^2\right)}\\
\sqrt{\!-\!4\left(\lambda^2\!-\!\frac{1}{4}(\alpha_4\!-\!\alpha_1)^2\right)
\left(\lambda^2\!-\!\frac{1}{4}(\alpha_4\!+\!\alpha_1)^2\right)}&
\!-\!2\left(\lambda\!+\!\frac{1}{2}(\alpha_4\!-\!\alpha_1)\right)
\left(\lambda\!-\!\frac{1}{2}(\alpha_4\!+\!\alpha_1)\right)
\end{bmatrix},\nonumber
\end{multline}
\begin{multline}
P_2=\cfrac{X\!-\!Y\!+\!Z}{2\alpha_2}+\frac{1}{2}I = \\
=\cfrac{1}{4\alpha_2\lambda}
\begin{bmatrix}
\!-\!2\lambda^2\!+\!2\alpha_2\lambda\!+\!\frac{1}{2}(\alpha_3^2\!-\!\alpha_2^2)&
e^{i\chi}\sqrt{\!-\!4\lambda^4\!+\!2(\alpha_2^2\!+\!\alpha_3^2)\lambda^2\!-\!
\frac{1}{4}(\alpha_3^2\!-\!\alpha_2^2)^2}\\
e^{\!-\!i\chi}\sqrt{\!-\!4\lambda^4\!+\!2(\alpha_2^2\!+\!\alpha_3^2)\lambda^2\!-\!
\frac{1}{4}(\alpha_3^2\!-\!\alpha_2^2)^2}&
2\lambda^2\!+\!2\alpha_2\lambda\!-\!\frac{1}{2}(\alpha_3^2\!-\!\alpha_2^2)
\end{bmatrix}=\\
=\cfrac{1}{4\alpha_2\lambda}
\begin{bmatrix}
\!-\!2\left(\lambda\!-\!\frac{1}{2}(\alpha_3\!+\!\alpha_2)\right)
\left(\lambda\!+\!\frac{1}{2}(\alpha_3\!-\!\alpha_2)\right)&
e^{i\chi}\sqrt{\!-\!4\left(\lambda^2\!-\!\frac{1}{4}(\alpha_3\!+\!\alpha_2)^2\right)
\left(\lambda^2\!-\!\frac{1}{4}(\alpha_3\!-\!\alpha_2)^2\right)}\\
e^{\!-\!i\chi}\sqrt{\!-\!4\left(\lambda^2\!-\!\frac{1}{4}(\alpha_3\!+\!\alpha_2)^2\right)
\left(\lambda^2\!-\!\frac{1}{4}(\alpha_3\!-\!\alpha_2)^2\right)}&
2\left(\lambda\!+\!\frac{1}{2}(\alpha_3\!+\!\alpha_2)\right)
\left(\lambda\!-\!\frac{1}{2}(\alpha_3\!-\!\alpha_2)\right)
\end{bmatrix},\nonumber
\end{multline}
\begin{multline}
P_3=\cfrac{X\!+\!Y\!-\!Z}{2\alpha_3}+\frac{1}{2}I = \\
=\cfrac{1}{4\alpha_3\lambda}
\begin{bmatrix}
\!-\!2\lambda^2\!+\!2\alpha_3\lambda\!-\!\frac{1}{2}(\alpha_3^2\!-\!\alpha_2^2)&
\!-\!e^{\!-\!i\chi}\sqrt{\!-\!4\lambda^4\!+\!2(\alpha_2^2\!+\!\alpha_3^2)\lambda^2\!-\!
\frac{1}{4}(\alpha_3^2\!-\!\alpha_2^2)^2}\\
\!-\!e^{\!-\!i\chi}\sqrt{\!-\!4\lambda^4\!+\!2(\alpha_2^2\!+\!\alpha_3^2)\lambda^2\!-\!
\frac{1}{4}(\alpha_3^2\!-\!\alpha_2^2)^2}&
2\lambda^2\!+\!2\alpha_3\lambda\!+\!\frac{1}{2}(\alpha_3^2\!-\!\alpha_2^2)
\end{bmatrix}=\\
=\cfrac{1}{4\alpha_3\lambda}
\begin{bmatrix}
\!-\!2\left(\lambda\!-\!\frac{1}{2}(\alpha_3\!+\!\alpha_2)\right)
\left(\lambda\!-\!\frac{1}{2}(\alpha_3\!-\!\alpha_2)\right)&
\!-\!e^{\!-\!i\chi}\sqrt{\!-\!4\left(\lambda^2\!-\!\frac{1}{4}(\alpha_3\!+\!\alpha_2)^2\right)
\left(\lambda^2\!-\!\frac{1}{4}(\alpha_3\!-\!\alpha_2)^2\right)}\\
\!-\!e^{\!-\!i\chi}\sqrt{\!-\!4\left(\lambda^2\!-\!\frac{1}{4}(\alpha_3\!+\!\alpha_2\right)^2)
\left(\lambda^2\!-\!\frac{1}{4}(\alpha_3\!-\!\alpha_2)^2\right)}&
2\left(\lambda\!+\!\frac{1}{2}(\alpha_3\!+\!\alpha_2)\right)
\left(\lambda\!+\!\frac{1}{2}(\alpha_3\!-\!\alpha_2)\right)
\end{bmatrix},\nonumber
\end{multline}
\begin{multline}
P_4=\cfrac{\!-\!X\!-\!Y\!-\!Z}{2\alpha_4}+\frac{1}{2}I = \\
=\cfrac{1}{4\alpha_4\lambda}
\begin{bmatrix}
2\lambda^2\!+\!2\alpha_4\lambda\!+\!\frac{1}{2}(\alpha_4^2\!-\!\alpha_1^2)&
\!-\!\sqrt{\!-\!4\lambda^4\!+\!2(\alpha_4^2\!+\!\alpha_1^2)\lambda^2\!-\!
\frac{1}{4}(\alpha_4^2\!-\!\alpha_1^2)^2}\\
\!-\!\sqrt{\!-\!4\lambda^4\!+\!2(\alpha_4^2\!+\!\alpha_1^2)\lambda^2\!-\!
\frac{1}{4}(\alpha_4^2\!-\!\alpha_1^2)^2}&
\!-\!2\lambda^2\!+\!2\alpha_4\lambda\!-\!\frac{1}{2}(\alpha_4^2\!-\!\alpha_1^2)
\end{bmatrix}=\\
=\cfrac{1}{4\alpha_4\lambda}
\begin{bmatrix}
2\left(\lambda\!+\!\frac{1}{2}(\alpha_4\!+\!\alpha_1)\right)
\left(\lambda\!+\!\frac{1}{2}(\alpha_4\!-\!\alpha_1)\right)&
\!-\!\sqrt{\!-\!4\left(\lambda^2\!-\!\frac{1}{4}(\alpha_4\!+\!\alpha_1)^2\right)
\left(\lambda^2\!-\!\frac{1}{4}(\alpha_4\!-\!\alpha_1)^2\right)}\\
\!-\!\sqrt{\!-\!4\left(\lambda^2\!-\!\frac{1}{4}(\alpha_4\!+\!\alpha_1)^2\right)
\left(\lambda^2\!-\!\frac{1}{4}(\alpha_4\!-\!\alpha_1)^2\right)}&
\!-\!2\left(\lambda\!-\!\frac{1}{2}(\alpha_4\!+\!\alpha_1)\right)
\left(\lambda\!-\!\frac{1}{2}(\alpha_4\!-\!\alpha_1)\right)
\end{bmatrix}.\nonumber
\end{multline}

\bigskip

{\bf 5.} Let us construct, according to p.~1, isometries
$\Gamma_{0,i}: H_i\to H_0,\; i=\ov{1,4}$ by projectors $P_i: H_0\to H_0$
and obtain regular indecomposable locally scalar representations of the graph
$\wt{D}_4$:
\begin{eqnarray*}
\Gamma_{0,1}=
\begin{bmatrix}
\sqrt{\cfrac{\lambda^2+\alpha_1\lambda-
(\alpha_4^2-\alpha_1^2)/4}{2\alpha_1\lambda}}\\
\sqrt{\cfrac{-\lambda^2+\alpha_1\lambda+
(\alpha_4^2-\alpha_1^2)/4}{2\alpha_1\lambda}}
\end{bmatrix}, &
\Gamma_{0,2}=
\begin{bmatrix}
\sqrt{\cfrac{-\lambda^2+\alpha_2\lambda+
(\alpha_3^2-\alpha_2^2)/4}{2\alpha_2\lambda}}\\
e^{-i\chi}\sqrt{\cfrac{\lambda^2+\alpha_2\lambda-
(\alpha_3^2-\alpha_2^2)/4}{2\alpha_2\lambda}}
\end{bmatrix}, \\
\Gamma_{0,3}=
\begin{bmatrix}
\sqrt{\cfrac{-\lambda^2+\alpha_3\lambda-
(\alpha_3^2-\alpha_2^2)/4}{2\alpha_3\lambda}}\\
-e^{-i\chi}\sqrt{\cfrac{\lambda^2+\alpha_3\lambda+
(\alpha_3^2-\alpha_2^2)/4}{2\alpha_3\lambda}}
\end{bmatrix}, &
\Gamma_{0,4}=
\begin{bmatrix}
\sqrt{\cfrac{\lambda^2+\alpha_4\lambda+
(\alpha_4^2-\alpha_1^2)/4}{2\alpha_4\lambda}}\\
-\sqrt{\cfrac{-\lambda^2+\alpha_4\lambda-
(\alpha_4^2-\alpha_1^2)/4}{2\alpha_4\lambda}}
\end{bmatrix}.
\end{eqnarray*}
$$
0 < (\alpha_4-\alpha_1)/2 \leq \lambda \leq
\min\left((\alpha_2+\alpha_3)/2,\,(\alpha_1+\alpha_4)/2\right),\quad
-\pi < \chi \leq \pi.
$$
Thus, indecomposable regular (not necessarily normalized) locally scalar
representations of the graph $\wt{D}_4$ depend on $6$ real parameters
(on the $4$ of $5$ parameters
$\alpha'_1,\,\alpha'_2,\,\alpha'_3,\,\alpha'_4,\,\alpha'_0$,
connected by the relation $\alpha'_1+\alpha'_2+\alpha'_3+\alpha'_4=2\alpha'_0,\;
\alpha_i=\cfrac{\alpha'_i}{\alpha'_0}$, and parameters $\lambda$ and $\chi$).

\end{document}